\documentclass{article}[11pt]

\usepackage{amsmath}
\usepackage{amssymb}
\usepackage{amsthm}
\usepackage{enumerate}
\usepackage{graphicx}
\usepackage{float}
\usepackage{array}
\usepackage{booktabs}

\newtheorem{theorem}{Theorem}
\newtheorem{lemma}[theorem]{Lemma}
\newtheorem{prop}[theorem]{Proposition}

\newtheorem{conj}[theorem]{Conjecture}
\newtheorem{definition}{Definition}[section]
\newtheorem{example}{Example}[section]
\newtheorem{question}{Question}[section]

\newcommand{\ZZ}{\mathbb{Z}}

\newcommand{\dic}{\text{Dic}_4}
\newcommand{\thickhline}{\specialrule{1.5pt}{1pt}{1pt}}
\newcolumntype{"}{@{\hskip\tabcolsep\vrule width 1.5pt\hskip\tabcolsep}}
\makeatother

\usepackage[table]{xcolor}
\usepackage{subfig}

\newcommand\RED{\cellcolor{red}}
\newcommand\BLUE{\cellcolor{blue}}

\newcommand\YEL{\cellcolor{yellow}}

\allowdisplaybreaks

\title{\bf Zarankiewicz Numbers\\
and Bipartite Ramsey Numbers\\[5ex]
}

\author{Alex F.~Collins\\
\small Rochester Institute of Technology\\[-0.8ex]
\small School of Mathematical Sciences\\[-0.8ex]
\small Rochester, NY 14623\\[-0.8ex]
\small \texttt{weincoll@gmail.com}\\[3ex]\and
Alexander W.~N.~Riasanovsky\\
\small University of Pennsylvania\\[-0.8ex]
\small Department of Mathematics\\[-0.8ex]
\small Philadelphia, PA 19104\\[-0.8ex]
\small \texttt{alexneal@math.upenn.edu}\\[3ex]\and
John C. Wallace\\
\small Trinity College\\[-0.8ex]
\small Department of Mathematics\\[-0.8ex]
\small Hartford, CT 06106\\[-0.8ex]
\small \texttt{john.wallace@trincoll.edu}\\[3ex]\and
Stanis{\l}aw P.~Radziszowski\\
\small Rochester Institute of Technology\\[-0.8ex]
\small Department of Computer Science\\[-0.8ex]
\small Rochester, NY 14623\\[-0.8ex]
\small \texttt{spr@cs.rit.edu}\\[3ex]
}

\date{\today}

\begin{document}
\maketitle
\thispagestyle{empty}

\begin{abstract}
The Zarankiewicz number $z(b;s)$ is the maximum size of a subgraph of $K_{b,b}$
which does not contain $K_{s,s}$ as a subgraph. The two-color bipartite Ramsey
number $b(s,t)$ is the smallest integer $b$ such that any coloring of the edges of
$K_{b,b}$ with two colors contains a $K_{s,s}$ in the first color or a $K_{t,t}$ in
the second color.

In this work, we design and exploit a computational method for bounding and
computing Zarankiewicz numbers.
Using it, we obtain several new values and bounds on $z(b;s)$ for $3 \le s \le 6$.
Our approach and new knowledge about $z(b;s)$ permit us to improve some of the
results on bipartite Ramsey
numbers obtained by Goddard, Henning and Oellermann in 2000. In particular, we
compute the smallest previously unknown bipartite Ramsey number, $b(2,5)=17$.
Moreover, we prove that up to isomorphism there exists a unique $2$-coloring which
witnesses the lower bound $16<b(2,5)$. We also find tight bounds on $b(2,2,3)$,
$17 \le b(2,2,3) \le 18$, which currently is
the smallest open case for multicolor bipartite Ramsey numbers.
\end{abstract}

\medskip
\noindent
{\bf Keywords:} Zarankiewicz number, bipartite Ramsey number\\
{\bf AMS classification subjects:} 05C55, 05C35
\eject

\section{Introduction}

\subsection*{Graph notation}

If $G$ is a bipartite graph, with the bipartition of its vertices
$V(G)=L(G)\cup R(G)$, or simply $V=L\cup R$, we will denote it by
writing $G[L,R]$. Furthermore, when we wish to point only to the orders
$m$ and $n$ of the {\em left} and {\em right} parts of the vertex
set $V(G)$, $m=|L|, n=|R|$, we will use notation $G[m,n]$.
The parts $L$ and $R$ will be called {\em left vertices} and
{\em right vertices} of $G$, respectively.
If $H$ is a subgraph of $G$, and its bipartition is
$H[L',R']$, then we will consider only the cases when
$L' \subset L$ and $R' \subset R$.
For the remainder of this paper, all bipartite graphs
have a fixed bipartition.

This allows us to treat any bipartite graph $G[m,n]$
as $m \times n$ 0-1 matrix $M_G$, whose rows are labeled
by $L$, columns are labeled by $R$, and where 1's stand for the
corresponding edges between $L$ and $R$.
The (bipartite) reflection of $G$ is obtained by swapping
the left and right vertices of $G$, or equivalently by transposing
the corresponding 0-1 matrix. The bipartite complement $\overline G$
of a (bipartite) graph $G[L,R]$, has the same bipartition as $G$,
but its matrix representation is the binary complement of $M_G$.  

\subsection*{Zarankiewicz numbers}

The {\em Zarankiewicz number} $z(m,n;s,t)$ is defined to be the maximum
number of edges in any subgraph $G[m,n]$ of the complete bipartite
graph $K_{m,n}$, such that $G[m,n]$ does not contain $K_{s,t}$.
For the diagonal cases, we will use $z(m,n;s)$ and $z(n;s)$
to denote $z(m,n;s,s)$ and $z(n,n;s,s)$, respectively.

In 1951, Kazimierz Zarankiewicz \cite{Zar} asked what is the minimum
number of 1's in a 0-1 matrix of order $n \times n$, which
guarantees that it has a $2 \times 2$ minor of 1's. In the notation
introduced above, it asks for the value of $z(n,n;2,2)+1$.

General Zarankiewicz numbers $z(m,n;s,t)$
and related extremal graphs have been
studied by numerous authors, including
K\"{o}v\'{a}ri, S\'{o}s, and Tur\'{a}n \cite{KST},
Reiman \cite{Rei}, Irving \cite{Irv},
and Goddard, Henning, and Oellermann \cite{GHO}.
A nice compact summary of what is known
was presented by Bollob\'as \cite{Bol} in 1995.
Recently, F\"uredi and Simonovits \cite{FuSi} published
an extensive survey of relationships between
$z(m,n;s,t)$ and much studied Tur\'{a}n numbers
${\rm ex} (k,K_{s,t})$.

The results and methods used to compute or estimate $z(n;2)$ are
similar to those in the widely studied case of ${\rm ex} (n,C_4)$,
where one seeks the maximum number of edges in any $C_4$-free
$n$-vertex graph. Previous papers established the exact values
of $z(n;s)$ for all $n \le 21$ \cite{DDR}, and some recent as
of yet unpublished work by Afzaly and McKay pushed it further to
all $n \le 31$ \cite{AfMc}, see also Table 3 in the Appendix.
Early papers by Irving \cite{Irv} and Roman \cite{Rom}
presented some bounding methods and results for concrete
cases with $s > 2$. For more data for $3 \le s \le 6$ see
our Appendix. For detailed discussion of general bounds
and asymptotics, especially for $s=2$ and $s=3$,
see the work by F\"uredi and Simonovits \cite{FuSi}.

\bigskip
\subsection*{Bipartite Ramsey numbers}

The {\em bipartite Ramsey number} $b(s_1,\dots, s_k)$ is
the least positive integer $b$ such that any coloring of
the edges of the complete bipartite graph $K_{b,b}$ with
$k$ colors contains $K_{{s_i},{s_i}}$ in the $i$-th
color for some $i$, $1\leq i\leq k$.

If $s_i=s$ for all $i$, then we
will denote this number by $b_k(s)$. The study of bipartite
Ramsey numbers was initiated by Beineke and Schwenk in 1976,
and continued by others, in particular Exoo \cite{Exoo},
Hattingh and Henning \cite{HaHe},
Goddard, Henning, and Oellermann \cite{GHO},
and Lazebnik and Mubayi \cite{LaMu}.

The connection between Zarankiewicz numbers and bipartite
Ramsey numbers is quite straightforward: the edges in color
$i$ in any coloring of $K_{n,n}$ witnessing
$n < b(s_1,\dots, s_k)$ give a lower bound witness
for $e \le z(n,n;s_i)$, where the $i$-th color has $e$ edges.
Thus, upper bounds on $z(n;s)$ can be useful in obtaining
upper bounds on bipartite Ramsey numbers.
This relationship was originally exploited by Irving \cite{Irv},
developed further by several authors, including
Goddard, Henning, and Oellermann \cite{GHO}, and it
will be used in this paper. The role of Zarankiewicz
numbers and witness graphs in the study of bipartite
Ramsey numbers is very similar to that of
Tur\'{a}n numbers ${\rm ex} (n,G)$ and $G$-free graphs
in the study of classical Ramsey numbers, where we color
the edges of $K_n$ while avoiding $G$ in some color.

For multicolor bipartite cases ($k>2$), we know most
when avoiding $C_4$, i.e. for $s=2$. The following
exact values have been established:
$b_2(2)=5$ \cite{BeSc}, $b_3(2)=11$ \cite{Exoo},
and $b_4(2)=19$ \cite{SP,DDR}. In the smallest
open case for 5 colors it is known that
$26 \le b_5(2) \le 28$ \cite{DDR}, where the lower
bound was obtained by a 5-coloring of $GF(5^2)\times GF(5^2)$,
and the upper bound is implied by a general upper bound
on $z(k^2+k-2)$ for $k=5$. It was also conjectured
that $b_5(2)=28$ \cite{DDR}.

Finally, we wish to point to the work by
Fenner, Gasarch, Glover and Purewal \cite{FGGP}, who wrote a very
extensive survey of the area of grid colorings, which are equivalent
to edge colorings of complete bipartite graphs. Their focus is
on the cases avoiding $C_4$ for both Zarankiewicz and Ramsey
problems.

\bigskip
\subsection*{Notes on asymptotics}

Asymptotics of Zarankiewicz numbers is quite well
understood (relative to Ramsey numbers).
The classical bound by
K\"{o}v\'{a}ri, S\'{o}s, and Tur\'{a}n \cite{KST},
generalized by several authors (cf. \cite{Bol,FuSi}) is
$$z(m,n;s,t)<
(s-1)^{1/t}(n-t+1)m^{1-1/t}+(t-1)m,
$$
which for constant $s=t$ becomes $z(n;s)=O(n^{2-1/t})$.
F\"uredi \cite{Fur} improved the general bound to
the best known so far
$$z(m,n;s,t)<
(s-t+1)^{1/t}(n-t+1)mn^{1-1/t}+(t-1)n^{2-2/t}+(t-2)m,
$$
for $m \ge t$ and $n \ge s \ge t \ge 2$.
These upper bounds are asymptotically optimal, as discussed
in a book chapter by Bollob\'as \cite{Bol}, and
a more recent monograph by F\"uredi and Simonovits \cite{FuSi}.

\medskip
In 2001, Caro and Rousseau \cite{CaRo}, using upper bounds
on Zarankiewicz numbers $z(n,n;s,s)$, proved that for
any fixed $m \ge 2$ there exist constants $A_m$ and $B_m$
such that for sufficiently large $n$ we have
$$A_m \Big({n \over {\log n}}\Big) ^{(m+1)/2} < b(m,n) <
B_m \Big({n \over {\log n}}\Big) ^m.$$
The asymptotics of other off-diagonal cases, including
avoidance of $K_{s,t}$ and other bipartite graphs,
was studied by Lin and Li \cite{LinLi}, and others.
For the diagonal case, the best known asymptotic
upper bound $b(n,n) < \big (1+o(1))2^{n+1} \log_2 n$ was
obtained by Conlon \cite{Con}.

\bigskip
\subsection*{Overview of this paper}

In the remainder of this paper we consider only the case of
avoiding  balanced complete $K_{s,s}$, i.e. the case of $s=t$.
Thus, for brevity, in the following the Zarankiewicz numbers
will be written as $z(m,n;s)$ or $z(n;s)$.

The main contribution of this paper is the method for
computing and bounding $z(m,n;s)$ for small $s > 2$,
and the results obtained by using it. The background
to the method and the method itself are presented in
Section 2.
The main results of this paper and the computations
leading to them are presented in Section 3.
The results are as follows: We obtain
several new values and bounds on $z(n;s)$ for $3 \le s \le 6$.
We compute the smallest previously unknown bipartite
Ramsey number, $b(2,5)=17$, and we prove that up to
isomorphism there exists a unique $2$-coloring which
witnesses the lower bound $16<b(2,5)$. Finally, we
find tight bounds on $b(2,2,3)$,
$17 \le b(2,2,3) \le 18$, which currently is
the smallest open case for multicolor bipartite Ramsey numbers.

\bigskip

\section{Two lemmas and their applications}
The focus of this section is on Lemmas \ref{star count} and \ref{density},
and refining their applications.  In the context of bipartite Ramsey
numbers and Zarankiewicz numbers, these lemmas may be found in various
forms in Section 12 of \cite{Guy}, in \cite{Irv}, \cite{GHO},
and \cite{DDR}. For use throughout the paper, we introduce
the following notation.  
\begin{definition}
Let $G[m,n]$ be some bipartite graph.  For positive integers $e,s,t$, we say that 
\begin{enumerate}
\item $G$ is a $(m,n,e^+)$-graph if $e(G)\geq e$,
\item $G$ is a $(m,n,e^+)_s$-graph if $G$ is a $(m,n,e^+)$-graph and $K_{s,s}\not\subseteq G$, and
\item $G$ is a $(m,n,e^+)_{s,t}$-graph if $G$ is a
$(m,n,e^+)_s$-graph and $K_{t,t}\not\subseteq \overline G$.
\end{enumerate}
In any of the above notations, we may
replace ``$e^+$'' with ``$e$'' whenever the condition
$e(G)\geq e$ is strengthened to $e(G) = e$ and drop ``$,e^+$'' whenever no restriction is placed on $e(G)$.  
\end{definition}
For example, ``$z(m,n;s)\geq z$'' is equivalent to ``there exists
some $(m,n,z)_s$-graph'' and ``$b(s,t)\geq m+1$'' is equivalent
to ``there exists some $(m,m)_{s,t}$-graph''. In general, we will
use the placeholder $\mathcal P$ to denote any of the empty word, ``$s$'', and ``$s,t$''.  
\begin{prop}\label{sum min}
For fixed positive integers $p,k,t$,
among the $k$-part sums $a_1+\cdots+a_k = p$ with $a_i\geq 0$, the sum
\begin{equation}\label{sum}
\sum_{i=1}^k{a_i\choose t}
\end{equation}
is minimized when $|a_i-a_j|\leq 1$ for all $1\leq i < j\leq k$.  
\end{prop}
\begin{proof}
Suppose $a_1 + \cdots + a_k = p$ is not balanced in the above sense,
and assume without loss of generality that $a_1\geq\cdots\geq a_k$.
It follows that $a_1-a_k>1$. Let $b_1 := a_1-1$,
$b_2 = a_2,\cdots,b_{k-1}=a_{k-1},b_k:=a_k+1$,
be a new split of $p$ into $k$ parts.
Note that
$\sum_{i=1}^k {a_i\choose t}-{b_i\choose t}=
{a_1-1\choose t-1}-{a_k\choose t-1}\geq0,$
thus a more balanced $k$-part sum does not increase (1).
Consequently, (1) is minimized for some sum as stated.
\end{proof}
\begin{lemma}[Star-Counting Lemma]\label{star count}
Let $G$ be a $(m, n, e^+)_s$-graph with
$e(G)=md_L+r_L=nd_R+r_R$, where $0\leq r_L<m$ and $0\leq r_R<n$.
If $(a_i)$ is the left degree sequence of $G$ $($and thus
$a_1 + \dots + a_m = e(G) )$, then 
\begin{equation}\label{star count 1}
(m-r_L)\cdot{d_L\choose s}+r_L\cdot{d_L+1\choose s}\leq\sum_{i=1}^m{a_i\choose s}\leq(s-1)\cdot{n\choose s}.  
\end{equation}
\end{lemma}
\begin{proof}
Fix $G$ as above. The first inequality in (\ref{star count 1}) follows
from Proposition \ref{sum min} since the leftmost expression corresponds
to a balanced $m$-part composition of $e$. Note that the middle expression
of (\ref{star count 1}) counts the number of stars $K_{1,s}\subseteq G$
whose center is on the left. If this sum exceeds $(s-1)\cdot{n\choose s}$,
then by the pigeonhole principle, there is some $B\subseteq R$, with
$|B|=s$ which is the set of leaves of at least $s$ of the left stars
as above. However, their union must contain a $K_{s,s}$ subgraph of $G$,
so the second inequality holds as well.  
\end{proof}
An analogous statement holds for the right-hand side of $G$.  
\begin{lemma}[Density Lemma]\label{density}
Let $G$ be $(m, n, e^+)_{\mathcal P}$-graph and $f=e-\lfloor e/m\rfloor$.
Then $G$ contains an induced $(m-1, n, f^+)_{\mathcal P}$-subgraph.
\end{lemma}
\begin{proof}
Fix $G$ as above and let $\overline{d}=e/m$.
Since $\overline d$ is the average left degree, we may find and remove some left vertex of degree at most $\overline{d}$, leaving us with a $(m-1, n, f^+)_{\mathcal{P}}$-graph.  
\end{proof}
An analogous statement holds for the right-hand side of $G$.  
Lemma 2 provides a static upper bound $z(m, n; s)\leq e$
based on the parameters $m, n, s, e$ alone.  
Likewise, given an upper bound on $z(m-1, n; s)$, Lemma
\ref{density} gives an upper bound on $z(m, n; s)$.  
\begin{example}
We will show that there is no $(4,4,10)_2$-graph.
The most balanced possible composition of $10$ into $4$ parts
is $2 + 2 + 3 + 3 = 10$. For $m = n = 4,\;s = 2$ and $d_L = r_L=2$
we have the left- and right-hand side of $(2)$ equal to $8$ and $6$,
respectively. Thus, $z(4;2)\leq 9$.
\end{example}

\begin{example}
We will show that there is no $(4,5,12+)_2$-graph.
If we apply Lemma $3$ to any $(5,4,12)_2$-graph, then we obtain
a $(4,4,10^+)$-graph. As argued in the previous example such
graphs do not exist, hence $z(4, 5; 2)\leq 11$.  
\end{example}
We combine Lemmas \ref{star count} and \ref{density} as follows.
\begin{lemma}\label{z comparison}
Let $m,n,s,w$ be positive integers.
If $z(m-1,n;s)<z\leq w-\lfloor w/m\rfloor$, then $z(m, n; s) < w$.
Also, if $a_1+\dots+a_m = w$ satisfies $|a_i-a_j|\leq 1$ for all
$1\leq i<j\leq m$ and ${a_1\choose s}+\cdots + {a_m\choose s}>(s-1){n\choose s}$,
then $z(m,n;s)<w$.  
\end{lemma}
Suppose we know that $W = (w_{ij})_{i,j=1}^{m,n}$ are upper bounds on
the Zarankiewicz numbers $z(i,j;s)$. Then, we may be able to improve
some of them using Lemma 4 by traversing the indices
$(i,j),\;1\leq i\leq m,\;1\leq j\leq n$ in some order
from $(1,1)$ to $(m,n)$ so that Lemma 4 can be applied
at each step. Call this algorithm \textsf{z\_bound}.  

\subsection*{Backwards paths extensions}\label{up down sec}
The bounds found just by \textsf{z\_bound} alone can be often improved
by exhaustive methods. If this is successful for any parameters,
further application of the \textsf{z\_bound} algorithm can
lead to improvements for higher parameters.
The same technique will be used to bound bipartite Ramsey numbers.
In both cases, we will attempt to construct all
$(m,n,e^+)_\mathcal P$-graphs (a possibly empty set).
To do this, we begin with all $(a,b,f^+)_\mathcal P$-graphs,
where $a,b,f$ are chosen carefully. For convenience, write
\[
(m,n,e^+)_{\mathcal P} \sqsupset (a,b,f^+)_{\mathcal P}
\]
if it is known that any $(m,n,e^+)_{\mathcal P}$-graph
contains some induced $(a,b,f^+)_{\mathcal P}$-subgraph.
When $(a,b)=(m-1,n)$ or $(a,b)=(m,n-1)$, the ``$\sqsupset$''
will be called a {\em step}.
A {\em backwards path} is a sequence of steps, such as
\[
(m_k, n_k, e_k^+) \sqsupset (m_{k-1}, n_{k-1}, e_{k-1}^+) \sqsupset \cdots \sqsupset (m_0, n_0, e_0^+).
\]
We aim at constructing the set of all $(m_k, n_k, e_k^+)_{\mathcal P}$-graphs,
up to bipartite graph isomorphism, using the following \textsf{extend} algorithm.
First, generate all of the $(m_0, n_0, e_0^+)_{\mathcal P}$-graphs
up to bipartite graph isomorphism by some other method.
Now suppose we have all $(m_i, n_i, e_i^+)_\mathcal P$-graphs.
For each such graph $G$, generate $(m_{i+1}, n_{i+1}, e_{i+1}^+)$-graphs
by adding a new vertex $v$ to the appropriate side of degree
$d(v)\geq e_{i+1}-e(G)$ in all possible ways. In addition, if
$\mathcal P = $``$s$'' or ``$s,t$'' and the bound
$z\geq z(m_{i+1},n_{i+1};s)$ is known, we may also
impose a condition that $d(v)\leq z-e(G)$.
Remove all generated graphs which are not of type
$(m_{i+1},n_{i+1},e_{i+1}^+)_\mathcal P$.
Reduce the remaining set up to bipartite
graph isomorphism, which can be readily accomplished by using
McKay's \textsf{nauty} package \cite{McK}.
Repeat this for all $G$ until all $(m_k, n_k, e_k^+)$-graphs are generated.  
\begin{question}\label{ques}
Given $a,b,m,n,e,\mathcal P$, how can we find a suitable backwards
path $(m, n, e^+)_\mathcal P \sqsupset \cdots \sqsupset (a, b, f^+)_\mathcal P$
in such a way that $f$ is as large as possible?
\end{question}
In principle, one could use techniques from dynamic programming to
obtain all such optimal paths, yet this is not practical in our case.
While Lemmas 2 and 3 are easy to apply at all times, the question
whether it is feasible to run \textsf{extend} algorithm depends
on other unpredictable factors.
\begin{example}
In order to aid in the computation of the bipartite Ramsey number
$b(2,5)$ and characterization of its lower bound witnesses, we found
a backwards path
$(16, 16, 189^+)_{5,2}\sqsupset\cdots\sqsupset (7, 7, f^+)_{5,2}$,
which is displayed in Figure $1$, where the rows and columns
correspond to $i$ and $j$ ranging from 7 to 16.
For reconstruction, the starting
$(7,7,42+)_{5,2}$-graphs have the number of edges close to
the maximum equal to $z(7;5)=44$ $($see Table $6)$.
\end{example}
\begin{figure}[H]
\centering
\includegraphics[scale=.75]{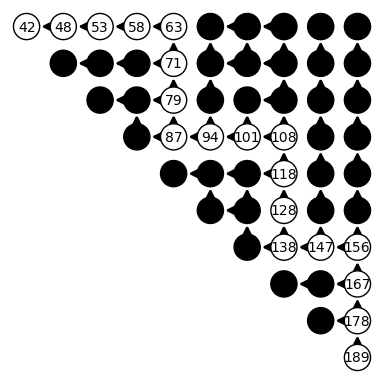}
\caption{Backwards path witnessing $(16,16,189^+)_{5,2} \sqsupset (7,7,42^+)_{5,2}$.}
\label{bt 16 16 189}
\end{figure}

The path highlighted in Figure \ref{bt 16 16 189} illustrates the sensitive
nature of this process. Parity plays a crucial role and it is not obvious
to the authors how in general to find a backwards path from $(m,n,e^+)$ to
the destination $(a,b,f^+)$ maximizing $f$ and feasible to follow with computations.
The pointers from entries indicate which immediately smaller parameters
where considered when performing computations leading to the displayed path.
If instead we step backwards along the main diagonal, we end at $(7,7,39^+)_{5,2}$.
Stepping back in two straight paths (straight to entry $(7,16)$, then straight
to entry $(7,7)$) coincidentally gives the same end value 39.  Note that a bare
density comparison gives only
$(16,16,189^+)_{5,2}\sqsupset (7,7,\left\lceil 49\cdot 189/256\right\rceil^+)_{5,2} = (7,7,37^+)_{5,2}$.
Up to isomorphism, there are
$7500$ $(7,7,37^+)_{5,2}$-graphs, $1619$ $(7, 7, 39^+)_{5,2}$-graphs, but only $33$ $(7, 7, 42^+)_{5,2}$-graphs.

\bigskip
\section{Bipartite Ramsey Numbers and Sidon Sets}\label{sids sec}

We motivate this section with two results from \cite{Col}. When trying
to establish the lower bound $16<b(2,5)$, one may consider searching for
witness graphs which satisfy certain structural properties.
The \textsf{nauty} package command 
$$
\textsf{genbg 16 16 64:64 -d4:4 -D4:4 -Z1}
$$
lists, up to isomorphism, all $4$-regular $(16,16,64)_2$-graphs where
the neighborhood of any two vertices intersect in at most one neighbor
(so that the graphs are $K_{2,2}$-free). This runs in a few minutes
on an ordinary laptop computer, and produces $19$ graphs. After
removing the graphs whose bipartite complement contains a $K_{5,5}$,
a single $(16,16,64)_{2,5}$-graph remains and its bipartite adjacency
matrix is shown in Table 1.

\medskip
\begin{table}[H]
\begin{center}
$\begin{array}{" c c c c " c c c c " c c c c " c c c c "}
\thickhline
  \RED& \RED& \RED& \RED& \BLUE& \BLUE& \BLUE& \BLUE& \BLUE& \BLUE& \BLUE& \BLUE& \BLUE& \BLUE& \BLUE& \BLUE  \\
    \BLUE& \BLUE& \BLUE&  \BLUE&  \RED& \RED& \RED& \RED& \BLUE& \BLUE&  \BLUE& \BLUE&  \BLUE& \BLUE&  \BLUE& \BLUE  \\
    \BLUE& \BLUE& \BLUE&  \BLUE& \BLUE&  \BLUE& \BLUE&  \BLUE&  \RED& \RED& \RED& \RED& \BLUE& \BLUE&  \BLUE&  \BLUE \\
    \BLUE& \BLUE& \BLUE&  \BLUE& \BLUE&  \BLUE& \BLUE&  \BLUE& \BLUE&  \BLUE& \BLUE&  \BLUE&  \RED& \RED& \RED& \RED\\\thickhline
  \RED& \BLUE& \BLUE&  \BLUE&  \RED& \BLUE& \BLUE&  \BLUE&  \RED& \BLUE& \BLUE&  \BLUE&  \RED& \BLUE& \BLUE& \BLUE  \\
    \BLUE& \RED& \BLUE& \BLUE&  \BLUE&  \RED& \BLUE& \BLUE&  \BLUE&  \RED& \BLUE& \BLUE&  \BLUE&  \RED& \BLUE& \BLUE  \\
    \BLUE& \BLUE&  \RED& \BLUE& \BLUE&  \BLUE&  \RED& \BLUE& \BLUE&  \BLUE&  \RED& \BLUE& \BLUE&  \BLUE&  \RED& \BLUE \\
    \BLUE& \BLUE& \BLUE&  \RED& \BLUE& \BLUE&  \BLUE&  \RED& \BLUE& \BLUE&  \BLUE&  \RED& \BLUE& \BLUE&  \BLUE&  \RED\\\thickhline
    \BLUE& \BLUE& \BLUE&  \RED& \BLUE& \BLUE&  \RED& \BLUE& \BLUE&  \RED& \BLUE& \BLUE&  \RED& \BLUE& \BLUE& \BLUE  \\
    \BLUE& \BLUE&  \RED& \BLUE& \BLUE&  \BLUE& \BLUE&  \RED& \RED& \BLUE& \BLUE&  \BLUE& \BLUE&  \RED& \BLUE& \BLUE  \\
    \BLUE& \RED& \BLUE& \BLUE&  \RED& \BLUE& \BLUE&  \BLUE& \BLUE&  \BLUE& \BLUE&  \RED& \BLUE& \BLUE&  \RED& \BLUE \\
  \RED& \BLUE& \BLUE&  \BLUE& \BLUE&  \RED& \BLUE& \BLUE&  \BLUE& \BLUE&  \RED& \BLUE& \BLUE&  \BLUE& \BLUE&  \RED\\\thickhline
    \BLUE& \RED& \BLUE& \BLUE&  \BLUE& \BLUE&  \BLUE&  \RED& \BLUE& \BLUE&  \RED& \BLUE&  \RED& \BLUE& \BLUE& \BLUE  \\
  \RED& \BLUE& \BLUE&  \BLUE& \BLUE&  \BLUE&  \RED& \BLUE& \BLUE&  \BLUE& \BLUE&  \RED& \BLUE&  \RED& \BLUE& \BLUE  \\
    \BLUE& \BLUE& \BLUE&  \RED& \BLUE&  \RED& \BLUE& \BLUE&  \RED& \BLUE& \BLUE&  \BLUE& \BLUE&  \BLUE&  \RED& \BLUE \\
    \BLUE& \BLUE&  \RED& \BLUE&  \RED& \BLUE& \BLUE&  \BLUE& \BLUE&  \RED& \BLUE& \BLUE&  \BLUE& \BLUE& \BLUE&    \RED\\\thickhline
\end{array}$
\caption{The bipartite adjacency matrix of a $16 < b(2,5)$ witness.  }\label{b25 Collins}
\end{center}
\end{table}
There is clear structure in this $(16,16,64)_{2,5}$-graph.
The matrix of Table 1 is a $4\times 4$ arrangement of
$4\times 4$ blocks, 12 of them being permutation matrices
of 4 elements. This graph has $2304=2^83^2$ automorphisms.

In \cite{Col}, a cyclic witness on 15 vertices to the $3$-color
bipartite Ramsey number $b(2,2,3)$ was found, but no witness of any
kind could  be found on 16 vertices.

\begin{question}\label{q Collins}
Are the bounds $16 < b(2,5)$ and $15 < b(2,2,3)$ tight?
\end{question}

In Theorem \ref{b25 17} of the next section, we will be able to
conclude that $b(2,5)=17$ and that the
lower bound witness found in Table \ref{b25 Collins} is indeed
the unique witness. For the second part of Question 3.1 we will
be able to improve the lower bound by 1 using bipartite Cayley
graphs as described in the remaining part of this section.

Our definition of bipartite Cayley graphs generalizes the classical
cyclic constructions. Given $\Gamma$ a group and $S\subseteq\Gamma$ a set of
\emph{edge generators}, the {\em bipartite Cayley graph generated by
$S$} is $X(\Gamma, S)$, where $V(X(\Gamma, S))=\Gamma\times\{1,2\}$ and
for each $g\in\Gamma$ and $s\in S$, there is an edge between $(g,1)$ and
$(g\cdot s,2)$. We impose no restriction on the symmetry of $S$ and
accept the identity $1_\Gamma$ as a valid edge generator.
We can easily describe what causes $X(\Gamma, S)$ to avoid $K_{2,2}$
using the concept of Sidon sets.
\begin{definition}
Given a group $\Gamma$, a subset $S\subseteq\Gamma$ is {\em Sidon} if there are
no solutions in $S$ to 
\[
s_1s_2^{-1}s_3s_4^{-1} = 1_\Gamma,
\]
unless $s_i = s_{i+1}$ for some $i = 0, 1, 2, 3$, with indices taken modulo $4$.  
\end{definition}
Sidon sets were originally defined over the integers,
while the above is a well-known generalization to arbitrary groups.
For a more detailed discussion, see \cite{HR} and \cite{TV}.
\begin{prop}
$X(\Gamma, S)$ is $K_{2,2}$-free if and only if $S$ is Sidon.  
\end{prop}
\begin{proof}
Suppose $S\subseteq \Gamma$ has a solution
$s_1,s_2,s_3,s_4\in S$ to $s_1s_2^{-1}s_3s_4^{-1}=1_\Gamma$.
Then with $a:=1_\Gamma$, $b:=s_1$, $c:=s_1s_2^{-1}$, and
$d:=s_1s_2^{-1}s_3 = s_4$, note that $a\neq c$ and $b\neq d$.
It follows that $K_{2,2}\subseteq X(\Gamma,S)$.   Conversely,
suppose there is some $K_{2,2}\subseteq X(\Gamma,S)$ with left
vertices $a,c$ and right vertices $b,d$. Then setting
$s_1:=a^{-1}b$, $s_2:=c^{-1}b$, $s_3:=c^{-1}d$ and $s_4:=a^{-1}d$,
we see that $s_1,s_2,s_3,s_4\in S$ satisfies
$s_1s_2^{-1}s_3s_4^{-1}=1_\Gamma$.
Assuming $a\neq c$ and $b\neq d$ (because of $s_i\neq s_{i+1}$),
this is genuinely a $K_{2,2}$.  
\end{proof}

The $3$-color construction on 15 vertices in \cite{Col} witnessing
$15 < b(2,2,3)$ can be described in terms of Sidon sets as follows:
Let $\Gamma$ be the additive group modulo 15, $\ZZ_{15}$, and consider
three bipartite Cayley graphs $X(\ZZ_{15}, S_i)$, $0\leq i\leq 2$, where
$S_0 = \{0, 1, 3, 7\}, S_1 = \{2, 4, 12, 13\}$, and
$S_2 = \{5, 6, 8, 9, 10, 11, 14\}$.
One can check that $S_0$ and $S_1$ are Sidon, $S_2$ yields
a $K_{3,3}$-free graph, and that the edges of $X(\ZZ_{15}, S_i)$'s
partition the edges of $K_{15,15}$.

We searched for witnesses to $16 < b(2,5)$ and $16 < b(2,2,3)$
using the same principle.
Our approach was to search among the bipartite Cayley graphs whose edge
generators are Sidon sets and whose groups are of order $16$.
This task is not difficult, since \textsf{sage} \cite{sage} has a
page\footnote{http://doc.sagemath.org/html/en/constructions/groups.html\#construction-instructions-for-every-group-of-order-less-than-32} which lists groups of small order alongside commands to generate them.
A few lines of code in \textsf{sage} allowed us to automate the process
of generating the Sidon sets $S$ in a group. Checking whether the bipartite
complement $X(\Gamma,\Gamma\setminus S)$ contains $K_{5,5}$ can also be
done easily in \textsf{sage}. We may assume that $1_\Gamma\in S$ without
loss of generality. Among the $14$ groups of order $16$, the only group
which produced a desired construction was $\dic$, the \emph{dicyclic}
group of order $16$ (a generalization of the quaternion group).
Up to isomorphism, only a single bipartite Cayley graph of the
form $X(\Gamma,S)$ witnessing $16<b(2,5)$ was found.
Interestingly, this graph is isomorphic to the one in
Table \ref{b25 Collins}.
Based on the same Sidon set, a 3-colored adjacency matrix of
$K_{16,16}$ corresponding to the bound $16<b(2,2,3)$ is
presented in Table \ref{b223 sage}.

\bigskip
\begin{table}[H]
\begin{center}
$\begin{array}{" c c c c " c c c c " c c c c " c c c c "}
\thickhline
\RED & \RED & \BLUE & \RED & \YEL & \YEL & \YEL & \YEL & \RED & \YEL & \BLUE & \BLUE & \YEL & \BLUE & \YEL & \YEL\\
\YEL & \RED & \RED & \BLUE & \RED & \YEL & \YEL & \YEL & \YEL & \RED & \YEL & \BLUE & \BLUE & \YEL & \BLUE & \YEL\\
\YEL & \YEL & \RED & \RED & \BLUE & \RED & \YEL & \YEL & \YEL & \YEL & \RED & \YEL & \BLUE & \BLUE & \YEL & \BLUE\\
\YEL & \YEL & \YEL & \RED & \RED & \BLUE & \RED & \YEL & \BLUE & \YEL & \YEL & \RED & \YEL & \BLUE & \BLUE & \YEL\\\thickhline
\YEL & \YEL & \YEL & \YEL & \RED & \RED & \BLUE & \RED & \YEL & \BLUE & \YEL & \YEL & \RED & \YEL & \BLUE & \BLUE\\
\RED & \YEL & \YEL & \YEL & \YEL & \RED & \RED & \BLUE & \BLUE & \YEL & \BLUE & \YEL & \YEL & \RED & \YEL & \BLUE\\
\BLUE & \RED & \YEL & \YEL & \YEL & \YEL & \RED & \RED & \BLUE & \BLUE & \YEL & \BLUE & \YEL & \YEL & \RED & \YEL\\
\RED & \BLUE & \RED & \YEL & \YEL & \YEL & \YEL & \RED & \YEL & \BLUE & \BLUE & \YEL & \BLUE & \YEL & \YEL & \RED\\\thickhline
\YEL & \BLUE & \BLUE & \YEL & \RED & \YEL & \YEL & \BLUE & \RED & \YEL & \YEL & \YEL & \YEL & \RED & \BLUE & \RED\\
\BLUE & \YEL & \BLUE & \BLUE & \YEL & \RED & \YEL & \YEL & \RED & \RED & \YEL & \YEL & \YEL & \YEL & \RED & \BLUE\\
\YEL & \BLUE & \YEL & \BLUE & \BLUE & \YEL & \RED & \YEL & \BLUE & \RED & \RED & \YEL & \YEL & \YEL & \YEL & \RED\\
\YEL & \YEL & \BLUE & \YEL & \BLUE & \BLUE & \YEL & \RED & \RED & \BLUE & \RED & \RED & \YEL & \YEL & \YEL & \YEL\\\thickhline
\RED & \YEL & \YEL & \BLUE & \YEL & \BLUE & \BLUE & \YEL & \YEL & \RED & \BLUE & \RED & \RED & \YEL & \YEL & \YEL\\
\YEL & \RED & \YEL & \YEL & \BLUE & \YEL & \BLUE & \BLUE & \YEL & \YEL & \RED & \BLUE & \RED & \RED & \YEL & \YEL\\
\BLUE & \YEL & \RED & \YEL & \YEL & \BLUE & \YEL & \BLUE & \YEL & \YEL & \YEL & \RED & \BLUE & \RED & \RED & \YEL\\
\BLUE & \BLUE & \YEL & \RED & \YEL & \YEL & \BLUE & \YEL & \YEL & \YEL & \YEL & \YEL & \RED & \BLUE & \RED & \RED\\
\thickhline
\end{array}$
\caption{A $3$-color bipartite adjacency matrix witnessing $16 < b(2,2,3)$.  }\label{b223 sage}
\end{center}
\end{table}

\section{Main Results}
\begin{theorem}\label{b25 17}
We have that
\[
b(2,5)=17,
\]
and the unique witness to $16<b(2,5)$ is the $(16,16)_{2,5}$-graph
given by Table $1$. Moreover, the only way to realize this witness
as a bipartite Cayley graph is with the dicyclic group $\dic$.  
\end{theorem}
\begin{proof}
The lower bound is implied by the constructions discussed in
previous section. The conclusion $b(2,5)=17$ follows by inspection
of Tables \ref{z 2} and \ref{z 5} in the Appendix:
$z(17;2)=74$ and $z(17;5)\leq 213$, adding to $287$, which
is not sufficient to cover all $289$ edges of $K_{17,17}$.
The uniqueness of the constructed graph on 16 vertices
follows from a similar but more detailed argument.
Since $z(16;2)=67$, it suffices to consider all of the
$(16,16,189^+)_{5,2}$-graphs. All of such graphs were
generated using the algorithms described in Section 2
along the computational backwards path displayed in
Figure 1.
The final computation (after performing many auxiliary computations
and consistency verifications) terminates in under a half hour on
an ordinary laptop computer. It returns a single graph
$G[16,16]$ with $192$ edges, isomorphic to the graph
given by Table \ref{b25 Collins}.
We wish to note that a seemingly much simpler approach
of considering all potential $(16,16,64^+)_{2,5}$-graphs
resulted to be computationally infeasible using our methods.
\end{proof}
\begin{theorem}\label{b223 bound}
It holds that
\[
17\leq b(2,2,3)\leq 18.
\]
\end{theorem}
\begin{proof}
The lower bound witness is found in Table \ref{b223 sage}.
The upper bound is implied by using the bounds in
Tables \ref{z 2} and \ref{z 3} in the Appendix:
$z(18;2)=81,\;z(18;3)\leq 156$, and $2\cdot 81 + 156 = 318 < 324$.  
\end{proof}
\begin{conj}
We conjecture that
$b(2,2,3)=17.$
\end{conj}
At present, Tables \ref{z 2} and \ref{z 3} are quite close to
providing a proof: $z(17;2)=74$ and $z(17;3)\leq 141$, so we have
that $2\cdot 74 + 141 =289$ is just barely too large.
It would suffice to prove that $z(17;3)\leq 140$,
but our computational attempts
to obtain this bound have proven to be too time-consuming.
The interested reader may note other weak-looking bounds
in Table \ref{z 3}, such as for $z(k,17;3)$
for $13\leq k\leq 17$.  

\bigskip
\section*{Acknowledgment}

We would like to thank the National Science Foundation Research Experiences
for Undergraduates Program (REU grant \#1358583) for support of the REU Site project,
which was held at the Rochester Institute of Technology during the summers
of 2013 and 2015. Most of the work reported in this paper was completed
during these REU summer sessions.

\bigskip
\bigskip
\section*{Appendix: Small Zarankiewicz Numbers}

The problem of computing $z(m,n;2)$ is well-studied (cf. \cite{Guy},
\cite{DDR}, \cite{AfMc}). Below in Table 3, we only list the values
of $z(n;2)$ until the first open case at $n=32$. More details on
$z(m,n;2)$ and related cases can be found in a recent work by
Afzaly and McKay \cite{AfMc}.

\begin{table}[H]
	\begin{center}
		$\begin{array}{c|c||c|c||c|c}
			n & z(n;2) & n & z(n;2) & n & z(n;2)\\
			\hline
			1 & 1 & 12 & 45 & 23 & 115\\
			2 & 3 & 13 & 52 & 24 & 122\\
			3 & 6 & 14 & 56 & 25 & 130\\
			4 & 9 & 15 & 61 & 26 & 138\\
			5 & 12 & 16 & 67 & 27 & 147\\
			6 & 16 & 17 & 74 & 28 & 156\\
			7 & 21 & 18 & 81 & 29 & 165\\
			8 & 24 & 19 & 88 & 30 & 175\\
			9 & 29 & 20 & 96 & 31 & 186\\
			10 & 34 & 21 & 105 & 32 & 189/190\\
			11 & 39 & 22 & 108 & 33\\
			\hline
		\end{array}$
		\caption{Zarankiewicz numbers $z(n;2)$ from \cite{DDR} and \cite{AfMc}.}\label{z 2}
	\end{center}
\end{table}

The following are tables of upper bounds on some small Zarankiewicz numbers.
A boldfaced entry is an exact value. A superscript $\,^\ast$ indicates
that there exists a unique $(m,n,z(m,n;s))_s$-graph. A superscript $\,^\dagger$
indicates that there is also a unique $(m,n,z(m,n;s)-1)_s$-graph. An italicized
entry indicates that the bound or value was determined with exhaustive computations.
Otherwise, an undecorated number indicates that the bound was obtained
by using Lemmas 2, 3 and 4, and without exhaustive enumeration
of $(m,n,e^+)_s$-graphs.

\bigskip
\begin{table}[H]
  \begin{center}
   $\begin{array}{c || *{13}{c|}}
        & 6 & 7 & 8 & 9 & 10 & 11 & 12 & 13 & 14 & 15 & 16 & 17 & 18\\\hline\hline
      6 & \mathbf{26}^\ast & \mathbf{29} & \mathbf{32} & \mathbf{36}^\ast & \mathbf{39}^\ast & \mathbf{42} & \mathbf{45}^\ast & \mathbf{48}^\ast & \mathbf{50} & \mathbf{53} & \mathbf{56} & \mathbf{58} & 61\\\hline
      7 & & \mathbf{33}^\ast & \mathbf{37}^\ast & \mathbf{40} & \mathbf{44}^\ast & \mathbf{47} & \mathbf{50} & \mathbf{53} & \mathbf{56} & \mathbf{60}^\ast & \mathbf{63}^\ast & \mathbf{66} & 69 \\\cline{0-0}\cline{3-14}
      8 & \multicolumn{2}{|c|}{} & \mathbf{42}^\ast & \mathbf{45} & \mathbf{50}^\ast & \mathbf{53} & \mathbf{57}^\ast & \mathbf{60} & \mathbf{64}^\ast & \mathbf{67} & \mathbf{70} & \mathbf{74}^\ast & 78 \\\cline{0-0}\cline{4-14}
      9 & \multicolumn{3}{|c|}{} & \mathbf{49} & \mathbf{54} & \mathbf{59}^\ast & \mathbf{64}^\ast & \mathbf{67}^\ast & \mathbf{70} & \mathbf{73} & \mathbf{77} & \mathbf{81} & 85      
\\\cline{0-0}\cline{5-14}
      10 & \multicolumn{4}{|c|}{} & \mathbf{60}^\dagger & \mathbf{64}^\ast & \mathbf{68} & \mathbf{73}^\ast & \mathbf{77} & \mathbf{81}^\ast & \mathbf{85}^\ast & \mathbf{90}^\ast & 94
      
\\\cline{0-0}\cline{6-14}
      11 & \multicolumn{5}{|c|}{} & \mathbf{69}^\ast & \mathbf{74} & \mathbf{80} & \mathbf{84} & \mathbf{88} & \mathbf{92} & \mathbf{96} & 101
\\\cline{0-0}\cline{7-14}
      12 & \multicolumn{6}{|c|}{} & \mathbf{80} & \mathbf{86}^\ast & \mathbf{91}^\ast & \mathbf{96} & \mathbf{99} & \mathbf{103}^\ast & 109 \\\cline{0-0}\cline{8-14}
      13 & \multicolumn{7}{|c|}{} & \mathbf{92}^\ast & \mathbf{98}^\ast & \mathbf{104}^\ast & \mathbf{107} & \mathit{110} & 116
\\\cline{0-0}\cline{9-14}
      14 & \multicolumn{8}{|c|}{} & \mathbf{105}^\ast & \mathbf{112}^\ast & \mathbf{115}^\ast & 118 & 124 \\\cline{0-0}\cline{10-14}
      15 & \multicolumn{9}{|c|}{} & \mathbf{120}^\dagger & \mathbf{123}^\ast & 126 & 132
\\\cline{0-0}\cline{11-14} 
      16 & \multicolumn{10}{|c|}{} & \mathbf{128}^\ast & \it{133} & 140 
\\\cline{0-0}\cline{12-14}
      17 & \multicolumn{11}{|c|}{} & 141 & 148 \\\cline{0-0}\cline{13-14}
      18 & \multicolumn{12}{|c|}{} & 156 \\\cline{0-0}\cline{14-14}
    \end{array}$
    \caption{Bounds on Zarankiewicz numbers $z(m,n;3)$.}\label{z 3}
  \end{center}
\end{table}

\bigskip
\begin{table}[H]
  \begin{center}
   $\begin{array}{c || *{13}{c|}}
        & 6 & 7 & 8 & 9 & 10 & 11 & 12 & 13 & 14 & 15 & 16 & 17 & 18\\\hline\hline
      6 & \mathbf{31}^\ast & \mathbf{36}^\ast & \mathbf{39} & \mathbf{43} & \mathbf{47} & \mathbf{51} & \mathbf{55} & \mathbf{59} & \mathbf{63} & \mathbf{67} & \mathbf{71}^\ast & \mathbf{75}^\ast & \mathbf{78}\\\hline
      7 & & \mathbf{42}^\dagger & \mathbf{45} & \mathbf{49} & \mathbf{54} & \mathbf{58} & \mathbf{63} & \mathbf{68}^\ast & \mathbf{72} & \mathbf{77} & \mathbf{82}^\ast & \mathbf{87^\ast} & \mathbf{90}
\\\cline{0-0}\cline{3-14}
      8 & \multicolumn{2}{|c|}{} & \mathbf{51}^\ast & \mathbf{55} & \mathbf{60} & \mathbf{65} & \mathbf{70} & \mathbf{75} & \mathbf{80} & \mathbf{85} & \mathbf{90} & \mathbf{95}^\ast & \mathbf{99} \\\cline{0-0}\cline{4-14}
      9 & \multicolumn{3}{|c|}{} & \mathbf{61} & \mathbf{67} & \mathbf{72} & \mathbf{78}^\ast & \mathbf{84}^\ast & \mathbf{88} & \mathbf{94} & \mathbf{99} & \mathbf{104} & \mathit{109}      
\\\cline{0-0}\cline{5-14}
      10 & \multicolumn{4}{|c|}{} & \mathbf{74}^\ast & \mathbf{79} & \mathbf{86}^\ast & \mathbf{93}^\ast & \mathbf{97} & \mathbf{103} & \mathbf{109} & \mathbf{115} & \mathbf{120} \\\cline{0-0}\cline{6-14}
      11 & \multicolumn{5}{|c|}{} & \mathbf{86} & \mathbf{93}^\ast & \mathbf{100}^\ast & \mathbf{105} & \mathbf{111} & \mathit{117} & 124 & 131 
\\\cline{0-0}\cline{7-14}
      12 & \multicolumn{6}{|c|}{} & 101 & 109 & 114 & 121 & 127 & 134 & 141 \\\cline{0-0}\cline{8-14}
      13 & \multicolumn{7}{|c|}{} & 118 & 123 & 131 & 137 & 145 & 152
\\\cline{0-0}\cline{9-14}
      14 & \multicolumn{8}{|c|}{} & 132 & 141 & 147 & 156 & 163 \\\cline{0-0}\cline{10-14}
      15 & \multicolumn{9}{|c|}{} & 151 & 157 & 166 & 174
\\\cline{0-0}\cline{11-14} 
      16 & \multicolumn{10}{|c|}{} & 167 & 177 & 185
\\\cline{0-0}\cline{12-14}
      17 & \multicolumn{11}{|c|}{} & 188 & 196 \\\cline{0-0}\cline{13-14}
      18 & \multicolumn{12}{|c|}{} & 207 \\\cline{0-0}\cline{14-14}
    \end{array}$
    \caption{Bounds on Zarankiewicz numbers $z(m,n;4)$. }\label{z 4}
  \end{center}
\end{table}

\bigskip
\begin{table}[H]
  \begin{center}
    $\begin{array}{c || *{13}{c|}}
        & 6 & 7 & 8 & 9 & 10 & 11 & 12 & 13 & 14 & 15 & 16 & 17 & 18\\\hline\hline
      6 & \mathbf{33}^\ast & \mathbf{38}^\ast & \mathbf{43}^\ast & \mathbf{48}^\ast & \mathbf{52} & \mathbf{57} & \mathbf{62} & \mathbf{67}^\ast & \mathbf{72}^\ast & \mathbf{76} & \mathbf{81} & \mathbf{86} & 91 \\\hline
      7 & & \mathbf{44}^\ast & \mathbf{50}^\ast & \mathbf{56}^\ast & \mathbf{60} & \mathbf{66} & \mathbf{72} & \mathbf{78}^\ast & \mathbf{84}^\dagger & \mathbf{88} & \mathbf{92}^\ast & \mathbf{96} & 101
\\\cline{0-0}\cline{3-14}
      8 & \multicolumn{2}{|c|}{} & \mathbf{57}^\ast & \mathbf{64}^\ast & \mathbf{68} & \mathbf{74} & \mathbf{80} & \mathbf{86}^\ast & \mathbf{92}^\ast & \mathbf{97} & \mathbf{103} & \mathbf{109} & 115
\\\cline{0-0}\cline{4-14}
      9 & \multicolumn{3}{|c|}{} & \mathbf{72}^\dagger & \mathbf{76} & \mathbf{82} & \mathbf{88} & \mathbf{95} & \mathbf{101} & \mathbf{108} & \mathbf{114} & \mathbf{121} & 128 \\\cline{0-0}\cline{5-14}
      10 & \multicolumn{4}{|c|}{} & \mathbf{84}^\ast & \mathbf{90} & \mathbf{97} & \mathbf{104} & \mathbf{110} & \mathbf{117} & \mathbf{124} & \mathbf{131} & 138 \\\cline{0-0}\cline{6-14}
      11 & \multicolumn{5}{|c|}{} & \mathbf{98} & \mathbf{106} & \mathbf{113} & \mathbf{120} & \mathbf{127} & \mathbf{135}^\ast & \mathbf{142} & 150 \\\cline{0-0}\cline{7-14}
      12 & \multicolumn{6}{|c|}{} & \mathbf{114} & \mathbf{122} & \mathbf{130} & \mathbf{138} & \mathit{146} & \mathbf{154}^\ast & 163 \\\cline{0-0}\cline{8-14}
      13 & \multicolumn{7}{|c|}{} &  \mathbf{132}^\ast & \mathbf{140} & \mathbf{149}^\ast  & \mathit{156}  & \mathit{165} & 174 \\\cline{0-0}\cline{9-14}
      14 & \multicolumn{8}{|c|}{} & \mathbf{150}^\ast & \mathbf{160}^\ast & \mathit{168} & \mathit{177} & 187 \\\cline{0-0}\cline{10-14}
      15 & \multicolumn{9}{|c|}{} & \mathbf{171}^\ast & \mathit{180} & \mathit{189} & 200 \\\cline{0-0}\cline{11-14}
      16 & \multicolumn{10}{|c|}{} & \mathbf{192}^\dagger & \mathit{201} & 212 \\\cline{0-0}\cline{12-14}
      17 & \multicolumn{11}{|c|}{} & \mathit{213} & 225 \\\cline{0-0}\cline{13-14}
      18 & \multicolumn{12}{|c|}{} & 238 \\\cline{0-0}\cline{14-14}
    \end{array}$
    \caption{Bounds on Zarankiewicz numbers $z(m,n;5)$.  }\label{z 5}
  \end{center}
\end{table}

\begin{table}[H]
  \begin{center}
   $\begin{array}{c || *{13}{c|}}
        & 6 & 7 & 8 & 9 & 10 & 11 & 12 & 13 & 14 & 15 & 16 & 17 & 18\\\hline\hline
      6 & \mathbf{35}^\ast & \mathbf{40} & \mathbf{45} & \mathbf{50} & \mathbf{55} & \mathbf{60} & \mathbf{65} & \mathbf{70} & \mathbf{75} & \mathbf{80} & \mathbf{85} & \mathbf{90} & \mathbf{95}\\\hline
      7 & & \mathbf{46}^\ast & \mathbf{52}^\ast & \mathbf{58}^\ast & \mathbf{64}^\ast & \mathbf{70}^\ast & \mathbf{75} & \mathbf{81} & \mathbf{87} & \mathbf{93} & \mathbf{99}^\ast & \mathbf{105}^\ast & \mathbf{110} \\\cline{0-0}\cline{3-14}
      8 & \multicolumn{2}{|c|}{} & \mathbf{59}^\ast & \mathbf{66}^\ast & \mathbf{73}^\ast & \mathbf{80}^\ast & \mathbf{85} & \mathbf{92} & \mathbf{99} & \mathbf{106} & \mathbf{113}^\ast & \mathbf{120}^\ast & \mathbf{125} \\\cline{0-0}\cline{4-14}
      9 & \multicolumn{3}{|c|}{} & \mathbf{74}^\ast & \mathbf{82}^\ast & \mathbf{90}^\ast & \mathbf{95} & \mathbf{102} & \mathbf{109} & \mathbf{116} & \mathbf{123} & \mathbf{130} & \mathbf{137}      
\\\cline{0-0}\cline{5-14}
      10 & \multicolumn{4}{|c|}{} & \mathbf{95}^\ast & \mathbf{100}^\ast & \mathbf{105} & \mathbf{112} & \mathbf{120} & \mathbf{127} & \mathbf{135} & \mathbf{142} & \mathbf{150} \\\cline{0-0}\cline{6-14}
      11 & \multicolumn{5}{|c|}{} & \mathbf{110}^\ast & \mathbf{115} & \mathbf{122} & \mathbf{130} & \mathbf{138} & \mathbf{147} & \mathbf{155} & \mathbf{163}
\\\cline{0-0}\cline{7-14}
      12 & \multicolumn{6}{|c|}{} & \mathbf{125}^\ast & \mathbf{132} & \mathbf{141} & \mathbf{150}^\ast & \mathbf{158} & \mathbf{167} & \mathbf{176} \\\cline{0-0}\cline{8-14}
      13 & \multicolumn{7}{|c|}{} & \mathbf{142} & \mathbf{152}^\ast & \mathbf{161} & \mathbf{170} & \mathbf{180} & \mathbf{189}
\\\cline{0-0}\cline{9-14}
      14 & \multicolumn{8}{|c|}{} & \mathbf{162} & \mathbf{172} & \mathit{182} & \mathbf{192}^\ast & \mathit{202} \\\cline{0-0}\cline{10-14}
      15 & \multicolumn{9}{|c|}{} & \mathbf{184}^\ast & 195 & 205 & 216
\\\cline{0-0}\cline{11-14} 
      16 & \multicolumn{10}{|c|}{} & 208 & 218 & 230 
\\\cline{0-0}\cline{12-14}
      17 & \multicolumn{11}{|c|}{} & 231 & 244 \\\cline{0-0}\cline{13-14}
      18 & \multicolumn{12}{|c|}{} & 258 \\\cline{0-0}\cline{14-14}
    \end{array}$
    \caption{Bounds on Zarankiewicz numbers $z(m,n;6)$.}\label{z 6}
  \end{center}
\end{table}

\end{document}